\documentclass[12pt]{amsart}
\usepackage{amssymb,amsmath,amsthm, color}
\usepackage{graphicx}
\newtheorem{theorem}{Theorem}

\newtheorem{corollary}[theorem]{Corollary}

\theoremstyle{remark}

\newtheorem*{remark}{Remark}

\numberwithin{theorem}{section} \numberwithin{equation}{section}

\newcommand{\Q}{\mathbb{Q}}

\newcommand{\Z}{\mathbb{Z}}

\begin{document}

\title[The 1729 $K3$ surface]
{The 1729 $K3$ surface}
    \author{Ken Ono and Sarah Trebat-Leder}
    \address{Department of Mathematics and Computer Science, Emory University,
    Atlanta, Georgia 30322}
    \email{ono@mathcs.emory.edu}
    \email{strebat@emory.edu}

    \thanks {The authors are supported by the NSF.  The first author also thanks  the Asa Griggs Candler Fund for their generous support.  The authors would also like to thank Bartosz Naskr\c{e}cki for correcting an error in an earlier version of the paper, and Doug Ulmer for his many helpful comments. 
    } 

\begin{abstract}
We revisit the mathematics that Ramanujan developed in connection with the famous ``taxi-cab" number $1729$. 
A study of his writings reveals that he had been studying Euler's diophantine equation
$$
a^3+b^3=c^3+d^3.
$$
It turns out that Ramanujan's work anticipated deep structures and phenomena which have become fundamental objects in arithmetic geometry and number theory.  We find that
he discovered a $K3$ surface with Picard number $18$, one which can be used to obtain infinitely many cubic twists over $\Q$ with rank $\geq 2$.
\end{abstract}
\maketitle
\noindent
\section{Introduction}
\subsection{1729 and Ramanujan}
Srinivasa Ramanujan is said to have possessed an uncanny memory for idiosyncratic properties of numbers. J. E. Littlewood famously remarked that ``every number is a personal friend of Ramanujan." This opinion is supported by the famous story of $1729$, the so-called ``Hardy-Ramanujan" number 
after the famous anecdote of G. H. Hardy concerning a visit to the hospital to see
his collaborator Srinivasa Ramanujan: (page 12 of \cite{Hardy}):
\medskip

\noindent
{\it 
I remember once going to see him when he was ill at Putney. I had ridden in taxi-cab number 1729 and remarked that the number seemed to me rather a dull one, and that I hoped it was not an unfavorable omen. ``No," he replied, ``it is a very interesting number; it is the smallest number expressible as the sum of two cubes in two different ways."}

\medskip
\noindent
Indeed, $1729$ is the smallest natural number which is the sum of two positive cubes in two different ways. We have that
\begin{equation}\label{TwoCubes}
1729\ = \ 9^3 + 10^3 \ = \ 12^3 + 1^3.
\end{equation}
While this anecdote might give one the impression that Ramanujan came up with this amazing property of $1729$ on the spot, he actually had written it down before even coming to England. 

Ramanujan, like many others throughout history, had been studying Euler's diophantine equation \begin{equation} \label{Diophantine}
X^3+Y^3=Z^3+W^3.\end{equation} For example: 

\begin{enumerate}
\item In 1913, he published the identity $$(6A^2 - 4AB + 4B^2)^3 = (3A^2 + 5AB - 5B^2)^3 + (4A^2 - 4AB + 6B^2)^3 + (5A^2 - 5AB - 3B^2)^3$$ in \cite{Ramanujan1}, which after dividing both sides by $27$ gives $12^3 = (-1)^3 + 10^3 + 9^3$.  
\item In entry 20 (iii) of Chapter 18 of his second notebook (see \cite{Berndt}), he gave the identity \begin{eqnarray*}& &(M^7 - 3M^4(1 +  P) + M(3(1 + P)^2 - 1))^3 \\ & &+ (2M^6 - 3M^3(1 + 2P) + (1 + 3P + 3P^2))^3 \\ & &+ (M^6 - (1 + 3P + 3P^2))^3 = (M^7 - 3M^4 P + M(3P^2 - 1))^3\end{eqnarray*} and explicitly wrote down (\ref{TwoCubes}) as part of a list of examples.
\item In 1915, he asked for solutions to $X^3 + Y^3 + Z^3 = U^6$ and $X^3 + Y^3 + Z^3 = 1$ in \cite{Ramanujan2} and \cite{Ramanujan3}. He included (\ref{TwoCubes}) as examples he wanted the parametrized solutions to account for.
\item On pg. 341 of his lost notebook (see \cite{AndrewsBerndt}), he offered a remarkable method for finding two infinite families of solutions $x_i^3 + y_i^3 = z_i^3 \pm 1$, which can be found in the appendix and involves expanding rational functions at zero and infinity.  He included (\ref{TwoCubes}) as part of a list of examples. 
\end{enumerate}

While Ramanujan's story has made (\ref{TwoCubes}) famous, it seems to have been first noticed by  B. Fr\'{e}icle de Bessyin $1657$ during his correspondences with Wallis and Fermat (See \cite{Dickson} for more on the history of this problem).    Euler later completely parametrized rational solutions to (\ref{Diophantine}), and while (1) and (2) above are not general, in his third notebook (pg 387) Ramanujan offered a family of solutions equivalent to Euler's general solution: 

If $$\alpha^2 + \alpha \beta + \beta^2 = 3 \lambda \gamma^2,$$ then 
$$(\alpha + \lambda^2 \gamma)^3 + (\lambda \beta + \gamma)^3 = (\lambda \alpha + \gamma)^3 + (\beta + \lambda^2 \gamma)^3$$

Although several other formulations equivalent to Euler's general solution have been discovered, many consider Ramanujan's to be the simplest of all.

\begin{remark}
His families in $(4)$ give infinitely many near misses to Fermat's Last Theorem for exponent $3$. 
\end{remark}

\begin{remark}In our modern perspective, the equation $X^3 + Y^3 = Z^3 + W^3$ gives a cubic surface, which is in fact a rational elliptic surface.  
\end{remark}

\subsection{Elliptic Curves}

Another classical question is which rational numbers $d$ can be written as the sum of two rational cubes.  For example, if $d = 1$, then this is equivalent to the exponent 3 case of Fermat's last theorem and was proven by Euler. By Ramanujan's time, there were many values of $d$ for which this equation was known to have no solutions.  For example, in the 1800's, Sylvester conjectured that if $p \equiv 5 \pmod{18}$ and $q \equiv 11 \pmod{18}$ are primes, then the equation has no solutions when $d = p, 2p, 4p^2, 4q, q^2, 2q^2$.  Pepin proved these claims and more (See \cite{Dickson} for more on the history of this problem).  

We now know that the equation $E: X^3 + Y^3 = 1$ is an elliptic curve with $j$-invariant 0, and for cube-free $d$, $E_d: X^3 + Y^3 = d$ is the cubic twist of $E$ by $d$.  Since $E_d$ is torsion free for all integral $d > 2$, this classical question is really one about ranks of elliptic curves in families of cubic twists. 

In particular, parametrized solutions to $X^3 + Y^3 = Z^3 + W^3$ give us families of elliptic curves with rank at least two.  Ramanujan gave us such solutions in (1).  Let 
\begin{eqnarray}\label{Points}P_1 &=& (x_1(T), y_1(T)) = (6T^2 - 4 T + 4, -3T^2 -5T + 5), \; \; \text{ and }\\ P_2 &=& (x_2(T), y_2(T)) = (4T^2 - 4T + 6, 5T^2 - 5T - 3) \nonumber
\end{eqnarray}
and $$k(T) = 63(3T^2 - 3T + 1)(T^2 + T+ 1)(T^2 - 3T + 3).$$
Then $$x_1(T)^3 + y_1(T)^3 = k(T) = x_2(T)^3 + y_2(T)^3.$$

Using these along with standard techniques in algebraic geometry, we can show the following: 

\begin{theorem}\label{rank}
The elliptic curve $$E_{k(T)}/\Q(T): X^3 + Y^3 = k(T)$$ has rank $2$.
\end{theorem}

\begin{remark}
While proving Theorem~\ref{K3Thm}, we will also show that the rank of $E_{k(T)}$ over $\overline{\Q}(T)$ is $4$.  In fact, rank 4 happens already over $\Q(\sqrt{-3})(T)$. 
\end{remark}

\begin{corollary}\label{twists}
All but finitely many $t \in \Q$ have that 
$E_{k(t)}/\Q$ has rank $\geq 2$. 
\end{corollary}

In a certain sense, elliptic curves with rank at least two are very rare.  Most specialists believe, based on both numerical evidence and heuristic considerations, that $50\%$ of all elliptic curves have rank 0 and $50\%$ have rank 1, leaving all curves of higher rank to fall in the remaining $0\%$.  More precisely, the Birch and Swinnerton-Dyer conjecture implies that the rank of an elliptic curve is equal to the ordering of vanishing of the associated $L$-series at $s = 1$.  This implies that the rank is even or odd according to the sign of the functional equation of this series, which is referred to as the parity conjecture. It is believed that this sign is $+$ or $-$ with equal frequency, and it is expected that almost all curves have the smallest rank compatible with the sign of their functional equation. However, numerical data of Zagier and Kramarz \cite{Zagier} suggests that this might not be the case for the family $X^3 + Y^3 = d$.  Their numerical data suggests that a positive proportion have rank at least 2, and a positive proportion have rank at least 3.  For $n = 2, 3$, let $S_n(X)$ be the number of cube free $d$ with $|d| \leq X$ and rank of $E_d$ at least $n$.  Assuming the parity conjecture, Mai \cite{Mai} proved $S_2(X) \gg X^{2/3 - \epsilon}$.   Stewart and Top \cite{StewartTop} then unconditionally showed that $S_2(X) \gg X^{1/3}$ and $S_3(X) \gg X^{1/6}$.  Assuming the parity conjecture, the exponent on their rank 3 result can be improved to $1/3$.  

\begin{remark}
We are able to tie some of Stewart and Top's results using our polynomial $k(t)$. In particular, using Theorem 1 of \cite{StewartTop}, we are able to get that $S_2(X) \gg X^{1/3}$.  Additionally, assuming the parity conjecture and using its explicit form for cubic twists as given in Section 12 of \cite{StewartTop}, we can also show that $S_3(X) \gg X^{1/3}$.   
\end{remark}

\subsection{$K3$-Surfaces}

We can also view $X^3 + Y^3 = k(T)$ as an elliptic surface, which turns out to be a $K3$ surface.  $K3$ surfaces, which were defined by Andr\'{e} Weil in 1958, have become fundamental objects in string theory, moonshine, arithmetic geometry, and number theory.  

\begin{theorem}\label{K3Thm}
The smooth minimal surface associated with the equation 
$$X^3 + Y^3 = k(T)$$ is an elliptic $K3$-surface with Picard number $18$ over $\overline{\Q}$.  
\end{theorem}

\begin{remark}
In Stewart and Top's case, the elliptic $K3$-surface has Picard number $20$.  
\end{remark}

\begin{remark}
In \cite{Dolgachev}, Dolgachev, van Geemen, and Kondo associate $K3$ surfaces to nodal cubic surfaces (i.e. cubic surfaces with at worst nodes as singularities).  These $K3$ surfaces admit a natural elliptic fibration and all have an automorphism of order $3$. Our $K3$ surface is part of this family.
\end{remark}
 
We will give the definition of a $K3$ surface and its Picard number in Section~\ref{K3_background}.

In Section~\ref{section_rank_proof}, we will prove Theorem~\ref{rank} and Corollary~\ref{twists}.   In Section~\ref{section_K3}, we will use Theorem~\ref{rank} and results on elliptic $K3$ surfaces to prove Theorem~\ref{K3Thm}. 

\section{Proof of Theorem~\ref{rank} and its corollary}\label{section_rank_proof}
\subsection{Background}
In this section, we will show that $E_{k(T)}/\Q(T)$ has rank two.  To study this rank, we use a map, described in Proposition 1 of  \cite{StewartTop}, from the $\Q(T)$ points of $E_{k(T)}$ to the vector space of holomorphic differentials on the auxiliary curve $C/\Q: S^3 = k(T)$.  For each point $P = (x(T), y(T))$ in $E_{k(T)}(\Q(T))$, we define an element $\phi_P$ of $\text{Mor}_\Q(C, E)$, where $E/\Q$ is given by $X^3 + Y^3 = k(T)$ as in the introduction, by $$\phi_P(T, S) = \left(\frac{x(T)}{S}, \frac{y(T)}{S}\right).$$ Then the map $$\lambda: E_{k(T)}(\Q(T)) \to H^0(C, \Omega^1_{C/\Q})$$ is given by $\lambda(P) = \phi_P^\ast \omega_E$, where $\phi^\ast_P \omega_E$ denotes the pullback via $\phi_P$ of the invariant differential $\omega_E$.  Proposition 1 of \cite{StewartTop} states that $\lambda$ is a homomorphism with finite kernel.  Therefore, we just need to understand its image.  
\subsection{Rank at least two}
To prove Theorem~\ref{rank}, we first need to show that the rank is at least 2.  As Ramanujan has given us two points, we just need to check their images and see that the differentials are linearly independent.   We compute that 
$$\phi_{P_1}^\ast \omega_E =  \frac{5S(6T + 5)}{4(3T^2 - 2T + 2)^2} \; dT \;\; \text{ and } \phi_{P_2}^\ast \omega_E =\frac{5S(2T - 1)}{4(2T^2 - 2T + 3)^2} \; dT,$$
where $P_1$ and $P_2$ are as in (\ref{Points}).

\subsection{Rank at most two}\label{L-function}
To show that the rank is at most two, we reduce the curve modulo a prime $p$.  Over a function field $\mathbb{F}_p(T)$, it is known that the rank is at most the order of vanishing of the $L$-function, and so we can get an upper bound on the rank of the reduced curve by computing its $L$-function in Magma.  For the prime of good reduction $p = 17$, the $L$-function factors as $$(17T - 1)^2(17T + 1)^2(83521T^4 + 34T^2 + 1),$$ which shows that the rank over $\mathbb{F}_p(T)$ is at most $2$.  See \cite{Ulmer} for more details on computing ranks of elliptic curves over function fields in this manner. 

Since reducing mod $p$ cannot decrease the rank, this shows that the rank over $\Q(T)$ is at most $2$.  

\subsection{Studying the cubic twists}

If $E/\Q(T)$ is an elliptic curve which is not isomorphic over $\Q(T)$ to an elliptic curve defined over $\Q$, we can specialize $T$ to a rational number $t$.  A result of Silverman \cite{Silverman} gives that for all but finitely many rational numbers $t$, this specialization map $\phi_{t}: E(\Q(T)) \to E_{t}(\Q)$ is an injective homomorphism.  Therefore, all but finitely many values of $t$ yield elliptic curves $E_{k(t)}(\Q)$ with rank at least 2.

\section{A $K3$ surface}
\label{section_K3}
\subsection{Background on K3 surfaces}\label{K3_background}
A $K3$ surface is a smooth minimal complete surface that is regular and has trivial canonical bundle.  Some examples of $K3$ surfaces include intersections of three quadrics in $\mathbb{P}^5$, intersections of a quadric and a cubic in $\mathbb{P}^4$, and a non-singular degree 4 surface in $\mathbb{P}^3$.  Andr\'{e} Weil named them in honor of Kummer, K\"{a}hler, Kodaira, and the mountain K2.  

The N\'{e}ron-Severi group of a variety $X$ is the group of divisors modulo algebraic equivalence, i.e. $$\text{NS}(X) := \text{Pic}(X)/\text{Pic}^0(X).$$  It is a finitely generated abelian group, and its rank is the Picard number of $X$.  For $K3$ surfaces over characteristic zero fields, this Picard Number is always $\leq 20$.

\subsection{Proof of Theorem~\ref{K3Thm}}
The criterion given in \cite{Beukers-Stienstra} on pages 276-277 show that it's a $K3$ surface.

To complete the proof, it suffices to compute the Picard number.  Using Tate's algorithm, this surface has six bad fibers, each of type IV.  The Shioda-Tata formula for the Picard number $\rho$ then says that 
$$\rho = r + 2 + 6 \cdot 2 = r + 14,$$
with $r$ the $\overline{\Q}(T)$-rank of the elliptic curve defined by our equation.  We know that the $\Q(t)$-rank is 2.  Also, note that $E_{(k(T)}$ has CM over $\Q(\sqrt{-3})$, and so the action of the endomorphism ring on our two independent points gives a $\Z$-module of rank 4. 

By our computation of the $L$-function for $E_{k(T)}$ reduced mod $17$ in (\ref{L-function}), we see that the rank over $\overline{\mathbb{F}}_{p}$ is at most $4$, and hence the rank over $\overline{\Q}(T)$ is too.  

Therefore, the rank over $\overline{\Q}(T)$ is exactly 4, and hence the Picard number is 18. 

\newpage
\section*{Appendix}

This is the famous page from Ramanujan`s Lost Notebook on which one finds his representations of $1729$ as the sum of two cubes in two different ways.

\bigskip
\begin{center}
\includegraphics[height=170mm]{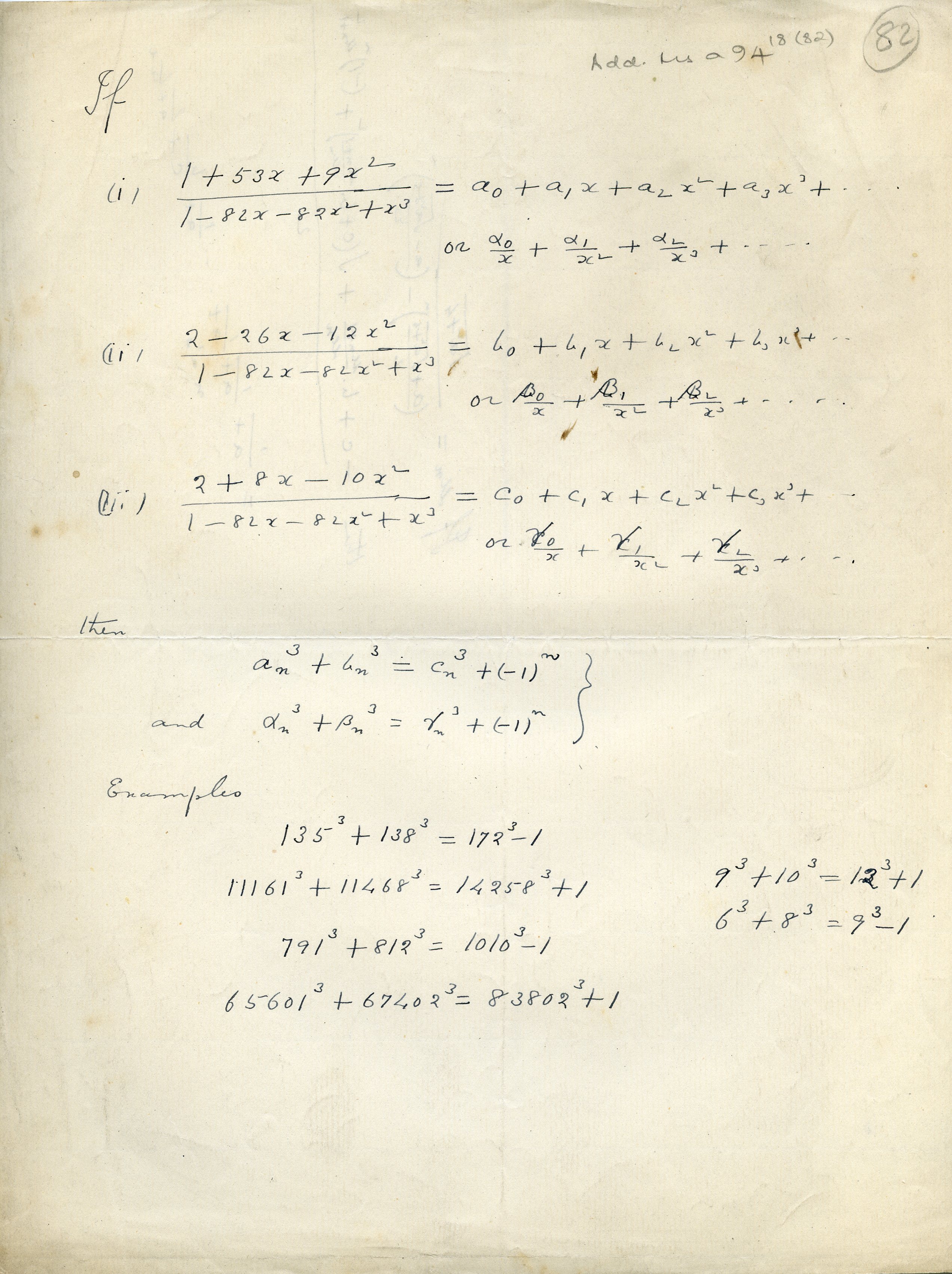}
\end{center}

\end{document}